\documentclass{amsart}

\makeatletter
 \def\LaTeX{\leavevmode L\raise.42ex
   \hbox{\kern-.3em\size{\sf@size}{0pt}\selectfont A}\kern-.15em\TeX}
\makeatother

\newcommand{\BibTeX}{{\rm B\kern-.05em{\sc
i\kern-.025emb}\kern-.08em\TeX}}
\newtheorem{theorem}{Theorem}[section]

\newtheorem{lemma}[theorem]{Lemma}
\theoremstyle{definition}
\newtheorem{definition}[theorem]{Definition}

\numberwithin{equation}{section}

\begin{document}

\title{Frames for spaces of Paley-Wiener functions on Riemannian manifolds}

\author{Isaac Pesenson}
\address{Department of Mathematics, Temple University,
Philadelphia, PA 19122} \email{pesenson@math.temple.edu}



 \begin{abstract}
It is shown that   Paley-Wiener functions on Riemannian manifolds
of bounded geometry can be reconstructed in a stable way from some
countable sets of  their inner products with certain distributions
of compact support. A reconstruction method in terms of frames is
given which is a generalization of the classical result of
Duffin-Schaeffer about exponential frames on intervals. All
results are specified in the case of the two-dimensional
hyperbolic space in its Poincare upper half-plane
 realization.

 \end{abstract}

\maketitle

 \section{Introduction}

A function $f\in L_{2}(R)$ is called $\omega$-bandlimited if its
$L_{2}$-Fourier transform
$$
\hat{f}(t)=\int_{-\infty}^{+\infty} f(x)e^{-2\pi i x t}dx
$$
has support in $[-\omega,\omega]$.

The Paley-Wiener theorem states that $f\in L_{2}(R)$ is
$\omega$-bandlimited if and only if
 $f$ is an entire function of exponential type not exceeding $2\pi
 \omega$. $\omega$-bandlimited functions form the Paley-Wiener class $PW_{\omega}$
  and often called  Paley-Wiener functions.

 The classical sampling theorem says, that if
$f$ is $\omega$-bandlimited then $f$ is completely determined by
its values at points $j/2\omega, j\in \mathbb{Z}$, and can be
reconstructed in a stable way from the samples $f(j/2\omega)$,
i.e.
 \begin{equation}
 f(x)= \sum_{j\in \mathbb{Z}} f\left(\frac{j}{2\omega}\right)\frac{\sin(2\pi\omega
(x-j/2\omega))}{2\pi \omega (x-j/2\omega )},
\end{equation}
where convergence is understood in the  $L_{2}$-sense. Moreover,
the following equality between "continuous" and "discrete" norms
holds true
\begin{equation}
\left(\int_{-\infty}^{+\infty}|f(x)|^{2}dt\right)^{1/2}=\left(\frac{1}{2\omega}\sum_{j\in
\mathbb{Z}}\left|f(j/2\omega)\right|^{2}\right)^{1/2}.
 \end{equation}

This equality follows from the fact that the functions $e^{2\pi i
t(j/2\omega)}$ form an orthonormal basis in
$L_{2}[-\omega,\omega]$.

The formulas (1.1) and (1.2) involve regularly spaced points
$j/2\omega, j\in \mathbb{Z}$. If one would like to consider
irregular sampling at a sequence of points $\{x_{j}\}$ and still
have a stable reconstruction from the samples $f(x_{j})$ then
instead of equality (1.2) the following Plancherel-Polya
inequality should hold true
\begin{equation}
C_{1}\sum_{j\in
\mathbb{Z}}|f(x_{j})|^{2}\leq\int_{-\infty}^{+\infty}|f(x)|^{2}dx\leq
C_{2}\sum_{j\in \mathbb{Z}}|f(x_{j})|^{2}.
\end{equation}
Such inequalities are also known as the frame inequalities.

Since the support of the Fourier transform $\hat{f}$ is in
$[-\omega, \omega]$  this inequality can be written in  the
following form
\begin{equation}
C_{1}\sum_{j\in\mathbb{Z}}\left|\int_{-\omega}^{+\omega}\hat{f}(t)e^{2\pi
i t
x_{j}}dt\right|^{2}\leq\int_{-\omega}^{+\omega}|\hat{f}(t)|^{2}dt\leq
C_{2}\sum_{j\in\mathbb{Z}}\left|\int_{-\omega}^{+\omega}\hat{f}(t)e^{2\pi
i t x_{j}}dt\right|^{2},
\end{equation}
that means that the functions $\left\{e^{2\pi i t x_{j}}\right\}$
form a kind of a basis (not necessary orthogonal), which is called
a \textit{frame} in the space $L_{2}[-\omega, \omega]$.

There  is a remarkable result of Duffin and Schaeffer \cite{DS},
that the inequalities (1.3) imply existence of a \textit{dual
frame} $\left\{\theta_{j}\right\}$  which consists of
$\omega$-bandlimited functions such that any $\omega$-bandlimited
function can be reconstructed according to the following formula
\begin{equation}
f(x)=\sum_{j\in \mathbb{Z}}f(x_{j})\theta_{j}(x),
\end{equation}
which is a generalization of the formula (1.1).

From this point of view the irregular sampling  was considered in
the classical paper of Duffin and Schaeffer \cite{DS} in which
they show that for the so-called uniformly dense sequences of
scalars $\{x_{j}\}, x_{j}\in \mathbb{R},$ the exponentials
$\{e^{ix_{j}t}\}$ form frames in appropriate spaces
$L_{2}[-\omega,\omega]$.
 In fact it was a far going
development of some ideas of Paley and Wiener \cite{PW34} about
irregular sampling. The theory of irregular sampling was very
active for many years \cite{B64}, \cite{BM67}, \cite{Lan67}, and
it is still active now \cite{OS}, \cite{LS02}.

The goal of the present article is to construct certain frames in
spaces of Paley-Wiener functions on Riemannian manifolds and to
show existence of reconstruction formulas of the type (1.5). Our
main Theorem 3.2 is a generalization of the Duffin-Schaeffer
result to the case of a Riemannin manifold of bounded geometry and
the Theorem 4.4 is a specification of this general result to the
case of the two-dimensional hyperbolic space in its upper-half
plane realization.

 The
notion of Paley-Wiener functions on manifolds were introduced in
the papers \cite{Pes89}-\cite{Pes04}. A subspace of Paley-Wiener
functions $PW_{\omega}(M)$
 on a Riemannian manifold  $M$ of
bounded geometry  consist   of all $L_{2}(M)$ functions whose
image in the spectral representation of the Laplace-Beltrami
operator $\Delta$ has support in the interval $[0,\omega]$.

In \cite{Pes00}, \cite{Pes04} a version of the Paley-Wiener
theorem was shown and an irregular sampling Theorem on manifolds
was proven. In the paper \cite{Pes01} a similar theory was
developed in the context of a general Hilbert space. The cases of
$\mathbb{R}^{d}$ and subelliptic cases on stratified Lie groups
were considered in \cite{Pes98a}, \cite{Pes98b}, \cite{Pes99}. In
all these situations the reconstruction algorithms were based one
the notion of variational splines on manifolds and in Hilbert
spaces. In the paper \cite{FP} an iterative reconstruction
algorithm was introduced.

The results of the present paper are different from our previous
results in the sense that

1) we use the notion of a frame to introduce a  method of
reconstruction of Paley-Wiener functions on manifolds;

2) we consider a kind of the "derivative sampling" which means
that we reconstruct functions from certain  sets of weighted
average values of $(1+\Delta)^{k} f, k\in \mathbb{N}.$

\bigskip

Note that interesting results in a similar direction on locally
compact groups were obtained recently by H. F\"uhr \cite{F} and H.
F\"uhr and K. Gr\"ochenig \cite {FG}.

 In the case of the hyperbolic space we  consider also reconstruction
 from pure "derivatives"  $\Delta^{k} f,$ for any fixed natural number
  $k>0$ . It is
 interesting to note that such result would be impossible to
 obtain on a general manifold. Moreover, this result does not hold
 true even in the case of $\mathbb{R}^{d}$. The reason why a
 Paley-Wiener function $f$ on the hyperbolic plane can be
 reconstructed from its pure derivative $\Delta^{k}f$ is, that
 the Laplace-Beltrami operator on the hyperbolic plane has bounded
 inverse. In this sense the situation on the hyperbolic plane is
 even "better" than on the $\mathbb{R}^{d}$. Note, that situation
 similar to the situation on the hyperbolic plane takes place on a
 general non-compact symmetric space.

Some partial results in these directions will appear in
\cite{Pes05}.

In the next section 2 some preliminary information on the subject
is given. The role of frames is explained in the section 3. An
example of the hyperbolic plane in its Poincare upper half-plane
realization is given in the section 4.

\section{Plancherel-Polya inequalities for Paley-Wiener functions
on Riemannian manifolds of bounded geometry}

 Let $M$, dim$M=d$, be
a connected $C^{\infty}-$smooth Riemannian manifold with a $(2,0)$
metric tensor $g$ that defines an inner product on every tangent
space $T_{x}(M), x\in M$. The corresponding Riemannian distance
$d$ on $M$is the function $d:M\times M\rightarrow R_{+}\bigcup
\{0\},$ which is defined as
$$
d(x,y)=inf\int_{a}^{b}
\sqrt{g(\frac{d\alpha}{dt},\frac{d\alpha}{dt})}dt,
$$
where $inf$ is taken over all $C^{1}-$curves
$\alpha:[a,b]\rightarrow M, \alpha(a)=x, \alpha(b)=y.$

 Let $ exp_ {x} $ : $T_{x}(M)\rightarrow M,$ be the exponential
geodesic map i. e. $exp_{x}(u)=\gamma (1), u\in T_{x}(M),$ where
$\gamma (t)$ is the geodesic starting at $x$ with the initial
vector $u$ : $\gamma (0)=x , \frac{d\gamma (0)}{dt}=u.$ If the
 $inj(M)>0  $ is the
injectivity radius of $M$ then
 the exponential map is a diffeomorphism of a ball of radius
$\rho<inj(M) $
 in the tangent space $T_{x}(M)$ onto the ball $B(x , \rho ).$ For
 every choice of an orthonormal (with respect to the inner
 product defined by $g$) basis of $T_{x}(M)$ the exponential map
  $exp$ defines a
 coordinate system on $B(x,\rho)$ which is called geodesic.

\bigskip

 Throughout the paper we will consider only geodesic coordinate
 systems.

\bigskip

We make the following assumptions about $M$:

\bigskip

1) the injectivity radius $inj(M)$ is positive;

\bigskip

2) for any $\rho\leq inj(M)$, and for every two canonical
coordinate systems $\vartheta_{x}: T_{x}(M)\rightarrow B(x,\rho),
\vartheta_{y}:T_{y}(M)\rightarrow B(x,\rho),$ the following holds
true
$$
\sup _{x\in B(x,\rho)\cap B(y,\rho)}\sup_{|\alpha|\leq k}|\partial
^{|\alpha|}\vartheta_{x}^{-1}\vartheta_{y}|\leq C(\rho, k);
$$
\bigskip

3) the Ricci curvature  $Ric$ satisfies (as a form) the inequality
\begin{equation}
Ric\geq-kg, k\geq 0.
\end{equation}

\bigskip

The Riemannian measure on $M$ is given in any coordinate system by
$$
d\mu=\sqrt{det(g_{ij})}dx,
$$
where the $\{g_{ij}\}$ are the components of the tensor $g$ in a
local coordinate system, and $dx$ is the Lebesgue's measure in
$R^{d}$.

The following Lemma was proved  in the paper \cite{Pes04}.

\begin{lemma}
For any Riemannian manifold of bounded geometry $M$  there exists
a natural $N_{M}$ such that  for any sufficiently small $\rho>0$
there exists a set of points $\{x_{i}\}$ such that
\begin{enumerate}
\item balls $B(x_{i}, \rho/4)$ are disjoint,

\item  balls $B(x_{i}, \rho/2)$ form a cover of $M$,

\item  multiplicity of the cover by balls $B(x_{i},\rho)$ is not
greater $N_{M}.$
\end{enumerate}
\end{lemma}

\bigskip

A set $\{x_{i}\}$ with such properties will be called a
$\rho$-lattice and will be denoted $Z(x_{i},\rho, N_{M})$.

 Let $K\subset B(x_{0},\rho/2)$ be a compact
subset and $\mu$ be a positive measure on $K$. We will always
assume that the total measure of $K$ is finite, i.e.
$$
0<|K|=\int_{K}d\mu<\infty.
$$

We consider the following distribution on
$C_{0}^{\infty}(B(x_{0},\rho)),$
\begin{equation}
\Phi(\varphi)=\int_{K}\varphi d\mu,
\end{equation}
where $\varphi \in C_{0}^{\infty}(B(x_{0},\rho)).$ As a compactly
supported distribution of order zero it has a unique continuous
extension to the space $C^{\infty}(B(x_{0},\rho))$.

Some examples of such  distributions which are of particular
interest to us are the following.

\bigskip

1) Delta functionals. In this case $K=\{x\}, x\in
B(x_{0},\rho/2),$ measure $d\mu$ is any positive number $\mu$ and
$\Phi(f)=\mu \delta_{x}(f)=\mu f(x).$

\bigskip

2) Finite or infinite sequences of delta functions $\delta_{j},
x_{j}\in B(x_{0}, \rho/2),$ with corresponding weights $\mu_{j}$.
In this case $K=\{x_{j}\}$ and
$$
\Phi(f)=\sum_{j}\mu_{j}\delta_{x_{j}}(f),
$$
where we assume the following
$$
0<|K|=\sum_{j}|\mu_{j}|<\infty, K=\{x_{j}\}.
$$

\bigskip

3) $K$ is a smooth submanifold in $B(x_{0},\rho/2)$ of any
codimension  and $d\mu$ is its "surface" measure.

\bigskip

4) $K$ is the closure of  $B(x_{0},\rho/2)$ and $d\mu$ is the
restriction of the Riemannian measure $dx$ on $M$.

\bigskip

We chose a lattice $Z(x_{i},\rho, N_{M})$ and in every ball
$B(x_{i},\rho/2)$ we consider a distribution $\Phi_{i}$ of type
(2.2) with support $K_{i}\subset B(x_{i}, \rho/2)$.

\bigskip

We say that a family $\Phi=\left\{\Phi_{j}\right\}$ is uniformly
bounded, if there exists a positive constant $C_{\Phi}$ such that
\begin{equation}
|K_{j}|\leq C_{\Phi}
 \end{equation}
 for all $j$.

We will also say that a family $\Phi=\left\{\Phi_{j}\right\}$ is
separated from zero if there exists a constant $c_{\Phi}>0$ such
that
\begin{equation}
|K_{j}|\geq c_{\Phi}
 \end{equation}
  for all $j$ where $|K_{j}|=\int_{K_{j}}d\mu_{j}$.

\bigskip

 To construct Sobolev spaces $W_{p}^{k}(M), 1\leq p\leq \infty,
k\in \mathbb{N},$ we consider a lattice $Z(y_{\nu}, \lambda,
N_{M}), \lambda<inj M$. We remained that  this assumption means in
particular that the multiplicity of the covers
$\{B(y_{\nu},\lambda)\}$ and $\{B(y_{\nu},\lambda/2)\}$ is not
greater $N_{M}$.

We construct a partition of unity  ${\varphi_{\nu}}$ that
subordinate to the family $\{B(y_{\nu},\lambda/2)\}$ and has the
following properties.

\bigskip

i)$\varphi_{\nu}\in C_{0}^{\infty} B(y_{\nu},\lambda/2 ),$

\bigskip

ii)$\sup_{x}\sup_{|\alpha|\leq k}|\varphi_{\nu}^{(\alpha
)}(x)|\leq C(k), $ where $ C(k) $ is independent on $\nu$ for
every $k $ in geodesic  coordinates.

\bigskip

We introduce Sobolev space $W_{p}^{k}(M), k\in \mathbb{N},$ as the
completion of $C_{0}^{\infty}(M)$ with respect to the norm
\begin{equation}
\|f\|_{W_{p}^{k}(M)}=\left(\sum_{\nu}\|\varphi_{\nu}f\|^{p}
_{W_{p}^{k}(B(y_{\nu},\lambda/2))}\right) ^{1/p}, k\in \mathbb{N}.
\end{equation}

We will use the notation $H^{k}(M)$ for the space $W_{2}^{k}(M)$.

Let $\Delta$ be the Laplace-Beltrami operator on a Riemannian
manifold $M, \dim M=d,$ with metric tensor $g$. In any local
coordinate system
$$
\Delta
f=\sum_{m,k}\frac{1}{\sqrt{det(g_{ij})}}\partial_{m}\left(\sqrt{det(g_{ij})}
g^{mk}\partial_{k}f\right)
$$
where $det(g_{ij})$ is the determinant of the matrix $(g_{ij})$.

 It is known that $\Delta$ is a self-adjoint positive
definite operator in the corresponding space $L_{2}(M,dx),$ where
$dx$ is the Riemannian measure. The regularity Theorem for the
Laplace-Beltrami operator $\Delta$ states that domains of the
powers
 $\Delta^{k/2}, k\in \mathbb{N},$ coincide with the Sobolev spaces
$H^{k}(M), k\in \mathbb{N},$ and  the norm (2.5) is equivalent to
the graph norm $\|f\|+\|\Delta^{k/2}f\|$.

\bigskip

We consider the positive square root $\Delta^{1/2}$ from the
positive definite self-adjoint operator $\Delta$. According to the
spectral theory \cite{BS} for a selfadjoint positive definite
operator $\Delta^{1/2}$ in a Hilbert space $L_{2}(M)$ there exist
a direct integral of Hilbert spaces $X=\int X(\lambda )dm (\lambda
)$ and a unitary operator $F$ from $L_{2}(M)$ onto $X$, which
transforms domain of $\Delta^{\mu/2}, \mu\geq 0,$ onto
$X_{\mu}=\{x\in X|\lambda ^{\mu}x\in X \}$ with norm

$$\|x(\lambda )\|_{X_{\mu}}= \left (\int^{\infty}_{0}
 \lambda ^{2\mu} \|x(\lambda )\|^{2}_{X(\lambda )} dm
 (\lambda ) \right )^{1/2}
 $$
 besides $F(\Delta^{\mu/2} f)=
 \lambda ^{\mu/2} (Ff), $ if $f$ belongs to the domain of
 $\Delta^{\mu/2}$.  As known, $X$ is the set of all $m $-measurable
  functions $\lambda \rightarrow x(\lambda )\in X(\lambda ) $, for which the norm

$$\|x\|_{X}=\left ( \int ^{\infty }_{0}\|x(\lambda )\|^{2}_{X(\lambda )}
 dm (\lambda ) \right)^{1/2} $$ is finite.

We will say that a function $f$ from $L_{2}(M)$ belong to the
Paley-Wiener space $PW_{\omega}$ if its "Fourier transform" $Ff$
has support in  $[0, \omega ] $.

The next two  theorems for an abstract selfadjoint operator in a
Hilbert space can be found in \cite{Pes00}, \cite{Pes01} .

\begin{theorem}
Let the $\mathcal{D}(\Delta^{k}), k\in \mathbb{N},$ be the domain
of the operator $\Delta^{k}$ and $\mathcal{D}^{\infty
}=\bigcap_{k\in \mathbb{N}}\mathcal{D}(\Delta^{k})$. The following
holds true:

\bigskip

  a)the set $\bigcup _{ \omega >0}PW_{\omega }(M)\subset
\mathcal{D}^{\infty }$ is dense in $L_{2}(M)$;

\bigskip

b) for every $\omega>0$ the set $PW_{\omega }(M)$ is a linear
closed subspace in $L_{2}(M)$.
\end{theorem}

Using the spectral resolution of identity  $P_{t}$ of the operator
$\Delta$ we define the unitary group of operators by the formula
$$
e^{it\Delta}f=\int_{0}^{\infty}e^{it\Delta}dP_{t}f, f\in L_{2}(M).
$$

 The next theorem can be considered
as a form of the Paley-Wiener theorem.
\begin{theorem}
The following conditions are equivalent:

1)$f\in PW_{\omega}(M)$;

2) for all $s\geq 0$ the following Bernstein inequality holds true
\begin{equation}
\|\Delta^{s}f\|\leq \omega^{2s}\|f\|;
\end{equation}

3) for every $g\in L_{2}(M)$ the scalar-valued function of the
real variable $t\in \mathbb{R}^{1}$
$$
\left<e^{it\Delta}f,g\right>=\int_{M}e^{it\Delta}f\overline{g}dx
$$
 is bounded on the real line and has an extension to the complex
plane as an entire function of the exponential type $\omega^{2}$;

4) the abstract-valued function $e^{it\Delta}f$  is bounded on the
real line and has an extension to the complex plane as an entire
function of the exponential type $\omega^{2}$.
\end{theorem}

The  following Theorem follows from the Sobolev embedding theorems
and from a well known result of Nelson \cite{N} about analytic
vectors.

\begin{theorem}
 The following continuous embeddings hold true
 \begin{equation}
 PW_{\omega}(M)\subset W_{p}^{m}(M),
 \end{equation}
 where
 $W_{p}^{m}(M)$, is the Sobolev space on $M$,
 $p>2,$ and in particular
 \begin{equation}  PW_{\omega}(M)\subset \mathrm{C}^{k}_{b}(M),
 \end{equation}
 where $C^{k}_{b}(M)$ is the space of $k$-differentiable bounded
 functions on $M$ and $k>d/2$.

 If, in addition, the manifold $M$ is real-analytic then every
 Paley-Wiener function is real-analytic.
\end{theorem}

 In our paper \cite{Pes04} the following
generalization of the Plancherel-Polya inequality was proved.

\begin{theorem}

For any given $ C_{\Phi}>0, c_{\Phi}>0, m=0,1,2,...,$ there exist
positive constants $C,  c_{1}, c_{2},$ such that for every
$\omega>0,$ every $\rho$-lattice $Z(x_{i},\rho, N_{M})$ with
$0<\rho< (C\omega)^{-1}$, every family of distributions
$\{\Phi_{i}\}$ of the form (2.2) with properties (2.3), (2.4) and
every $f\in PW_{\omega}(M)$ the following inequalities hold true

\begin{equation}
c_{1}\left(\sum_{j}\left|\Phi_{j}(f)\right|^{2}\right)^{1/2}\leq
\rho^{-d/2}\|f\|_{L_{2}(M)} \leq c_{2}\left(\sum_{j}
|\Phi_{j}(f)|^{2}\right)^{1/2}.
\end{equation}

\end{theorem}

\bigskip

 In the case of Euclidean
space when $\Phi_{i}=\delta_{x_{i}}$ and $\{x_{i}\}$ is the
regular lattice the above inequality represents the so-called
Plancherel-Polya inequality.

We prove the following extension of the Theorem 2.5.

\begin{theorem}

For any given $C_{\Phi}>0, c_{\Phi}>0, m=0,1,2,...,$ there exist
positive constants $C,  c_{1}, c_{2},$ such that for every
$\omega>0, $ every $\rho$-lattice $Z(x_{i},\rho, N_{M})$ with
$0<\rho< (C\omega)^{-1}$, every family of distributions
$\{\Phi_{i}\}$ of the form (2.2) with properties (2.3), (2.4) and
every $f\in PW_{\omega}(M)$ the following inequalities hold true

$$
c_{1}\left(\sum_{j}\left|\Phi_{j}\left((1+\Delta)^{n}f\right)\right|^{2}\right)^{1/2}\leq
\rho^{-d/2}\|f\|_{L_{2}(M)} \leq
$$
\begin{equation}
c_{2}\left(\sum_{j}
\left|\Phi_{j}\left((1+\Delta)^{n}f\right)\right|^{2}\right)^{1/2}.
\end{equation}

This implies in particular that on the space $PW_{\omega}(M)$ the
$L_{2}(M)$ norm is equivalent to the norm of the Sobolev space
$H^{m}(M)$ for every $m=0,1,2,...$.

\end{theorem}

\begin{proof}
Indeed, since the space $PW_{\omega}(M)$ is invariant under the
Laplace-Beltrami operator we have according to the Theorem 2.5
\begin{equation}
\rho^{-d/2}\|(1+\Delta)^{n}f\|_{L_{2}(M)} \leq
c_{2}(n)\left(\sum_{j}
\left|\Phi_{j}\left((1+\Delta)^{n}f\right)\right|^{2}\right)^{1/2}.
\end{equation}
But, because the operator $(1+\Delta)^{n}$ has bounded inverse
\begin{equation}
\|f\|_{L_{2}(M)}=\left\|(1+\Delta)^{-n}(1+\Delta)^{n}f\right\|_{L_{2}(M)}\leq
c(k)\left\|(1+\Delta)^{n}f\right\|_{L_{2}(M)}.
\end{equation}
These two inequalities imply the second part of the inequality
(2.10).

By the same inequality (2.9) we have
$$
c_{1}\left(\sum_{j}\left|\Phi_{j}\left((1+\Delta)^{n}f\right)\right|^{2}\right)^{1/2}\leq
\rho^{-d/2}\left\|(1+\Delta)^{n}f\right\|_{L_{2}(M)},
$$
and then the Bernstein inequality for $f$ gives the left side of
the inequality (2.10).
\end{proof}

\section{Uniqueness, stability and reconstruction in terms of  frames}

We consider the set of distributions of the form

$$
\Phi_{j}^{(n)}=\left(1+\Delta\right)^{n}\Phi_{j}, n\in
\mathbb{N}\bigcup\{0\},
$$
where $ n\in \mathbb{N}\bigcup\{0\},$ is a fixed number and
$$
 \left(1+\Delta\right)^{n}\Phi_{j}(f)=\Phi_{j}\left(\left(1+\Delta\right)^{n}f\right)=
 \int_{K_{j}}\left(1+\Delta\right)^{n}f d\mu_{j}.
 $$

 We say that a set of functionals $\Phi^{(n)}=\left\{\Phi_{j}^{(n)}
\right\}$ is a uniqueness set for $PW_{\omega}(M)$, if every $f\in
PW_{\omega}(M)$ is uniquely determined by its values
$\left\{\Phi_{j}^{(n)}(f)\right\}$.

For any such set $\Phi^{(n)}$ and any $\omega>0$ the notation
$l_{2}^{\omega}\left(\Phi^{(n)}\right)$ will be used for a linear
subspace of all sequences $\{v_{j}\}$ in $l_{2}$ for which there
exists a function $f$ in $PW_{\omega}(M)$ such that
$$
\Phi_{j}^{(n)}(f)=v_{j}.
$$
In general $l_{2}^{\omega}\left(\Phi^{(n)}\right)\neq l_{2}$.

\begin{definition}
A linear reconstruction method $R$ for a set
$\Phi^{(n)}=\left\{\Phi_{j}^{(n)}\right\}$ is a linear operator
$$
R:l_{2}^{\omega}\left(\Phi^{(n)}\right)\rightarrow
PW_{\omega}\left(M\right)
$$
such that
$$
R:\left \{\Phi_{j}^{(n)}(f)\right\}\rightarrow f.
$$

 The reconstruction method is said to be stable, if it is
continuous in the topologies induced respectively by $l_{2}$ and
$L_{2}(M)$.
\end{definition}

\begin{theorem}
For the given $n\in \mathbb{Z}, C_{\Phi}>0,c_{\Phi}>0$, there
exists a constant $C>0$ such that for any $\omega>0$, any
$\rho<\left(C\omega\right)^{-1}$,  any lattice $Z(x_{j},\rho,
N_{M})$, and any family of distributions $
\Phi=\left\{\Phi_{j}\right\}$ of the type (2.2) that satisfy (2.3)
and (2.4) with the given $C_{\Phi}, c_{\Phi}$, every function
$f\in PW_{\omega}(M)$ is uniquely defined by the set of samples
$\left\{\Phi_{j}\left((1+\Delta)^{n}f\right)\right\}_{j\in\mathbb{Z}}$,
in other words, if for a $\{v_{j}\}\in l_{2}$ there exists a
function $f\in PW_{\omega}(M)$ such that
$\Phi_{j}\left((1+\Delta)^{n}f\right)=v_{j}$ for all $j$, then
such function is unique.

Moreover, any reconstruction method from the set of samples
$\left\{\Phi_{j}^{(n)}(f)\right\}_{j\in\mathbb{Z}}$ is stable.
\end{theorem}

The proof of the Theorem is an immediate consequence of the
Plancherel-Polya inequalities (2.10).

Next, using the idea and the method of Duffin and Schaeffer
\cite{DS} we are going to describe  a stable method of
reconstruction of a function $f\in PW_{\omega}(M)$ from the
samples
$$
\left\{\Phi_{j}\left((1+\Delta)^{n}f\right)\right\}_{j\in\mathbb{Z}}\in
l_{2}
$$.

\begin{theorem}

For the given $n\in \mathbb{Z}, C_{\Phi}>0,c_{\Phi}>0$, there
exists a constant $C>0$ such that for any $\omega>0$, any
$\rho<\left(C\omega\right)^{-1}$,  any lattice $Z(x_{j},\rho,
N_{M})$, and any family of distributions $
\Phi=\left\{\Phi_{j}\right\}$ of the type (2.2) that satisfy (2.3)
and (2.4) with the given $C_{\Phi}, c_{\Phi}$, the following
statement holds true:

there exists a  frame
 $\{\Theta_{j}^{(n)}\}$ in the space $PW_{\omega}(M)$
such that every $\omega$-band limited function $f\in
PW_{\omega}(M)$ can be reconstructed from a set of samples
$\left\{\Phi_{j}\left((1+\Delta)^{(n)}(f)\right)\right\}\in l_{2}$
by using the formula
\begin{equation}
f=\sum_{j}\Phi_{j}\left((1+\Delta)^{(n)}(f)\right)\Theta_{j}^{(n)}.
\end{equation}

\end{theorem}

\begin{proof}
For the functional
$$
f\rightarrow \Phi_{j}\left((1+\Delta)^{(n)}(f)\right),
$$
defined on the space $PW_{\omega}(M)$ we will use the notation
$\Phi_{j}^{(n)}$. The definition of the functionals $\Phi_{j}$
(2.2) and the Bernstein inequality for functions from the space
$PW_{\omega}(M)$ imply that  every such functional is continuous
on $PW_{\omega}(M)$. By the Riesz theorem there are functions
$\phi_{j}^{(n)}\in PW_{\omega}(M)$ such, that for every $f\in
PW_{\omega}(M)$
$$
\Phi_{j}^{(n)}(f)=\left<\phi_{j}^{(n)},f\right>=\int_{M}\phi_{j}^{(n)}\overline{f}
$$
Moreover, the inequalities (2.10) show that the set of functionals
$\left\{\Phi_{j}^{(n)}\right\}$ is a frame in the space
$PW_{\omega}(M)$.

  The next goal is to show that the
so-called frame operator
\begin{equation}
Ff=\sum_{j}\left<\phi_{j}^{(n)},f\right>\phi_{j}^{(n)},\>\>\>\>
f\in PW_{\omega}(M),
\end{equation}
is an automorphism of the space $PW_{\omega}(M)$ onto itself and
$$
 \|F\|\leq c_{2},\>\>\>\> \|F^{-1}\|\leq c_{1}^{-1},
$$
where $c_{1}, c_{2}$ are from (2.10). Let us introduce the
operator
$$
F_{J}:PW_{\omega}(M)\rightarrow  PW_{\omega}(M),
$$
which is given
by the formula
$$
F_{J}f=\sum_{j\leq
J}\left<f,\phi_{j}^{(n)}\right>\phi_{j}^{(n)},\>\>\>\> f\in
PW_{\omega}(M).
$$
By the frame inequalities (2.10) and the Holder inequality we have
$$
\|F_{J_{1}}f-F_{J_{2}}f\|^{2}=\sup_{\|h\|=1}
\left|\sum_{J_{1}<j\leq
J_{2}}\left<f,\phi_{j}^{(n)}\right>\left<\phi_{j}^{(n)},h\right>\right|^{2}\leq
c_{2} \sum_{J_{1}<j\leq
J_{2}}\left|\left<f,\phi_{j}^{(n)}\right>\right|^{2}.
$$
By the same frame inequality the right side goes to zero when
$J_{1},J_{2}$ go to infinity. Thus the limit
$$
\lim _{J\rightarrow \infty}F_{J}f=Ff,\>\>\>\> f\in PW_{\omega}(M),
$$
does exists. Next,
$$
\|Ff\|^{2}=\sup_{\|h\|=1}
\left|\sum_{j}\left<f,\phi_{j}^{(n)}\right>\left<\phi_{j}^{(n)},h\right>
\right|^{2}\leq \sup_{\|h\|=1}c_{2}^{2}\|f\|^{2}\|h\|^{2}
=c_{2}^{2}\|f\|^{2},
$$
which shows that the operator $F$ is continuous.

 Now, the frame inequalities (2.10) imply that
$$
c_{1}I\leq F\leq c_{2}I,
$$
where $I$ is the identity operator. Thus, we have
$$
c_{2}^{-1}F\leq I,\>\>\> c_{2}^{-1}F^{-1}\geq c_{1}c_{2}^{-1}I,
$$
and then
$$
0\leq I-c_{2}^{-1}F\leq
I-c_{1}c_{2}^{-1}I=\left(c_{2}-c_{1}\right)c_{2}^{-1}I.
$$
It implies
$$
\left\|I-c_{2}^{-1}F\right\|\leq
\left\|\left(c_{2}-c_{1}\right)c_{2}^{-1}I\right\|\leq
\left(c_{2}-c_{1}\right)c_{2}^{-1}<1.
$$
 Thus, it shows that the operator
$\left(c_{2}^{-1}F\right)^{-1}$ and consequently the operator
$F^{-1}$ are well defined bounded operators and because
$$
F^{-1}=c_{2}^{-1}\left(C_{2}^{-1}F\right)^{-1}=c_{2}^{-1}
\sum_{m=0}^{\infty}\left(I-c_{2}^{-1}F\right)^{m},
$$
it gives us the desired estimate $\|F^{-1}\|\leq c_{1}^{-1}$.

We obtain
$$
f=F^{-1}Ff=F^{-1}\lim_{J\rightarrow \infty}\sum_{j\leq
J}\left<f,\phi_{j}^{(n)}\right>\phi_{j}^{(n)}=
\sum_{j}\left<f,\phi_{j}^{(n)}\right>\Theta_{j}^{(n)},
$$
where
$$
\Theta_{j}^{(n)}=F^{-1}\phi_{j}^{(n)},
$$
gives a dual frame $\Theta_{j}^{(n)}$ in the space
$PW_{\omega}(M)$.

 The Theorem is
proven.

\end{proof}

 \section{The Poincare upper half-plane}

As an illustration of our results we  will consider the hyperbolic
plane in its upper half-plane realization.

  Let $G=SL(2,\mathbb{R})$ be the special linear group of all
$2\times 2$ real matrices with determinant 1 and $K=SO(2)$ is the
group of all rotations of $\mathbb{R}^{2}$. The factor
$\mathbb{H}=G/K$ is known as the 2-dimensional hyperbolic space
and can be described in many different ways. In the present paper
we consider the realization of $\mathbb{H}$ which is called
Poincare upper half-plane (see \cite{Helg1}, \cite{H}, \cite{T}).

As a Riemannian manifold $\mathbb{H}$ is identified with the
regular upper half-plane of the complex plane
$$
\mathbb{H}=\{x+iy|x,y\in \mathbb{R}, y>0\}
$$
with a new  Riemannian metric
$$
ds^{2}=y^{-2}(dx^{2}+dy^{2})
$$
and corresponding Riemannian measure
 $$
  d\mu=y^{-2}dxdy.
   $$
   If we
define the action of $\sigma\in G$ on $z\in \mathbb{H}$ as a
fractional linear transformation
$$
\sigma\cdot z=(az+b)/(cz+d),
$$
then the metric $ds^{2}$ and the measure $d\mu$  are invariant
under the action of $G$ on $\mathbb{H}$. The point $i=\sqrt{-1}\in
\mathbb{H}$ is invariant for all $\sigma\in K$. The Haar measure
$dg$ on $G$ can be normalizes in a way that the  following
important formula holds true
$$
\int_{\mathbb{H}}f(z)y^{-2}dxdy=\int_{G}f(g\cdot i) d g.
$$

In the corresponding space of square integrable functions
$L_{2}(G)$ with the inner product
$$
<f,h>=\int_{\mathbb{H}}f\overline{h}y^{-2}dxdy
$$
we consider the Laplace-Beltrami operator

 $$
 \Delta=y^{2}\left(\partial_{x}^{2}+\partial_{y}^{2}\right)
 $$
of the metric $ds^{2}$.

 It is known that $\Delta$ as an operator
in $L_{2}(\mathbb{H})=L_{2}(\mathbb{H},d\mu)$ which initially
defined on $C_{0}^{\infty}(\mathbb{H})$ has a self-adjoint closure
in $L_{2}(\mathbb{H})$.

Moreover, if $f$ and $\Delta f$ belong to $L_{2}(\mathbb{H})$,
then
$$
<\Delta f,f>\leq -\frac{1}{4}\|f\|^{2},
$$
where $\|f\|$ means the $L_{2}(\mathbb{H})$ norm of $f$.

 The Helgason transform
of $f$ for $s\in \mathbb{C}, \varphi\in (0,2\pi]$, is defined by
the formula
$$
\hat{f}(s,\varphi)=\int_{\mathbb{H}}f(z)\overline{Im(k_{\varphi}z)^{s}}
y^{-2}dxdy,
$$
where $k_{\varphi}\in SO(2)$ is the rotation of $\mathbb{R}^{2}$
by angle $\varphi$.

We have the following inversion formula for all $f\in
C_{0}^{\infty}(\mathbb{H})$
$$
f(z)=(8\pi^{2})^{-1} \int_{t\in\mathbb{R}}\int_{0}^{2\pi}
\hat{f}(it+1/2,\varphi)Im(k_{\varphi}z)^{it+1/2}t \tanh \pi t
d\varphi dt.
$$

The Plancherel Theorem states that a such defined  map
$f\rightarrow \hat{f}$ can be extended to an isometry of
$L_{2}(\mathbb{H})$ with respect to invariant measure $d\mu $ onto
$L_{2}(\mathbb{R}\times (0,2\pi ])$ with respect to the measure
$$
\frac{1}{8\pi^{2}} t \tanh \pi t dt d\varphi.
$$
\bigskip

If $f$ is a function on $\mathbb{H}$ and $\varphi$ is a
$K=SO(2)$-invariant function on $\mathbb{H}$ their convolution is
defined by the formula
$$
f\ast \varphi(g\cdot i)=\int_{SL(2,\mathbb{R})}f(gu^{-1}\cdot
i)\varphi(u)du, i=\sqrt{-1},
$$
where $du$ is the Haar measure on $SL(2,\mathbb{R})$. It is known,
that for the Helgason-Fourier transform the following formula
holds true
$$
\widehat{f\ast\varphi}=\hat{f}\cdot\hat{\varphi}.
$$

The following formula holds true
\begin{equation}
\widehat{\Delta f}=-\left(t^{2}+\frac{1}{4}\right)\hat{f}.
\end{equation}
and our Theorem 2.3 takes the following form.
\begin{theorem}
A function $f\in L_{2}\left(\mathbb{H}\right)$ belongs to the
space $PW_{\omega}(\mathbb{H})$ if and only if for every
$\sigma\in \mathbb{R}$ the following holds true
\begin{equation} \|\Delta^{\sigma} f\|\leq
\left(\omega^{2}+\frac{1}{4}\right)^{\sigma}\|f\|.
 \end{equation}
 \end{theorem}
As it follows from the Theorem 2.4 if $f\in
PW_{\omega}(\mathbb{H})$ then for every $\sigma\geq 0$ the
function $\Delta^{\sigma} f$ belongs to $C^{\infty}(\mathbb{H})$
and is bounded on $\mathbb{H}$.

 Moreover  functions from $PW_{\omega}(\mathbb{H})$ are
not just infinitely  differentiable on
 $ \mathbb{H}$ but they are
\textit{real analytic functions } on the upper half-plane.

Indeed, every $f\in PW_{\omega}(\mathbb{H})$ is an analytic vector
of the elliptic differential operator $\Delta$ which means
\cite{N} that the series
$$
\sum_{k=0}^{\infty}\frac{\|\Delta^{k}f\|}{k!}s^{k}\leq
e^{\left(\omega^{2}+\frac{1}{4}\right) s}\|f\|
$$
is convergent for any $s>0$. By the famous result  of Nelson
\cite{N} it implies that the function $f$ is real analytic on
$\mathbb{H}$. It shows in particular that the support of the
function $f$ is the entire half-plane $\mathbb{H}$.

This fact can be treated  as a
 form of the uncertainty principle on $\mathbb{H}$:

\bigskip

\textit{If the support of the Helgason-Fourier transform
$\widehat{f}$ of a function $f\in L_{2}(\mathbb{H}, y^{-2}dxdy)$
is contained in a set $(-\omega, \omega)\times (0,2\pi]$, then the
support of $f$ is the entire half-plane $\mathbb{H}$}.

\bigskip

By applying the Theorem 2.3 we obtain the following result.
\begin{theorem}
A function $f$ is an $\omega$-band-limited signal, if for any
$g\in L_{2}(\mathbb{H},d\mu)$ the complex valued function
$$
t\rightarrow
<e^{it\Delta}f,d\mu>=\int_{\mathbb{H}}e^{it\Delta}f(z)
\overline{g(z)}y^{-2}dxdy
$$
of the real variable $t$ is an entire function of exponential type
$\omega^{2}+\frac{1}{4}$ which is bounded on the real line.
\end{theorem}

 It is known that the one parameter group of operators $e^{it\Delta }$ acts
on  functions by the following formula
$$
e^{it\Delta }f(z)=f\ast G_{t},
$$
where
$$
G_{t}(ke^{-r}i)=(4\pi)^{-1}\int_{\eta\in
\mathbb{R}}e^{-i(\eta^{2}+1/4)t}P_{i\eta-1/2}(cosh r)\eta \tanh
\pi\eta d\eta,
$$
here $k\in SO(2),$ $r$ is  the geodesic distance, $ke^{-r} i$ is
representation of points of $\mathbb{H}$ in the geodesic polar
coordinate system on $\mathbb{H}$, and $P_{i\eta -1/2}$ is the
associated Legendre function.

The last Theorem can be reformulated in the following terms.
\begin{theorem}
A function $f\in L_{2}(\mathbb{H}, d\mu)$ is $\omega$-band-limited
if and only if  for every $g\in L_{2}(\mathbb{H}, d\mu)$ function
$$
t\rightarrow <f\ast G_{t},g>=\int_{\mathbb{H}}f\ast
G_{t}\overline{g} d\mu
$$
is an entire function of the exponential type
$\omega^{2}+\frac{1}{4}$ which is bounded on the real line.
\end{theorem}

We are going to describe our reconstruction algorithm using the
language of frames.

 Note that according to the inversion
formula for the Helgason-Fourier transform we have

$$
\Phi_{j} f(z)=(8\pi^{2})^{-1} \int_{t\in\mathbb{R}}\int_{0}^{2\pi}
\hat{f}(it+1/2,\varphi)\Phi_{j}\left(Im(k_{\varphi}z_{j})^{it+1/2}\right)t
\tanh \pi t d\varphi dt.
$$
which means that the Helgason-Fourier transform of a distribution
 $\Phi_{j}$ is given by the formula
\begin{equation}
\widehat{\Phi_{j}} =
\Phi_{j}\left(Im(k_{\varphi}z)^{it+1/2}\right)
\end{equation}
and
$$
\widehat{\Delta^{n}\Phi_{j}} =
\left(\omega^{2}+\frac{1}{4}\right)^{n}\Phi_{j}\left(Im(k_{\varphi}z)^{it+1/2}\right)
$$

In the case of Poincare upper half-plane it is not difficult to
formulate our main Theorem 3.2. Moreover as the formula 4.1 shows
the Laplacian $\Delta$ in the space $L_{2}(\mathbb{H},
y^{-2}dxdy)$ has bounded inverse which allows reconstruction from
pure derivatives $\Delta^{n}f$. Namely
\begin{theorem}

Given two constants $0<c_{\Phi}\leq C_{\Phi}$ and an $n\in
\mathbb{N}\bigcup\{0\}$, there exists a constant $c=c(C_{\Phi},
c_{\Phi}, N, n)>0$ such that for any $\omega>0$ , any
$(\rho,N)$-lattice $Z\left(x_{j},\rho,N\right)$ with
$$
0<\rho<c\left(\omega^{2}+\frac{1}{4}\right)^{-1/2},
$$
and every family of distributions $\Phi=\{\Phi_{j}\}$ that satisfy
(2.2)-(2.4) with the given $c_{\Phi}, C_{\Phi}$, the following
statements hold true.

1)  The set of analytic functions
 $\{\widehat{\Delta^{n}\Phi_{j}}\}$ is
  a frame in the space
$$
L_{2}\left([-\omega, \omega]\times (0,2\pi ], \frac{1}{8\pi^{2}} t
\tanh \pi t dt d\varphi\right).
$$

 2) There exists a  frame
 $\{\Theta_{j}^{(n)}\}$ in the space $PW_{\omega}(\mathbb{H})$
such that every $\omega$-band limited function $f\in
PW_{\omega}(\mathbb{H})$ can be reconstructed from a set of
samples $\left\{\Phi_{j}(\Delta^{n}f)\right\}$ by using the
formula
$$
f=\sum_{j}\Phi_{j}\left(\Delta^{n}f\right)\Theta_{j}^{(n)}.
$$

\end{theorem}

When $n=0$ and every $\Phi_{j}$ is a Dirac measure
$\delta_{x_{j}}$ at a point $x_{j}\in \mathbb{H}$, the Theorem can
be considered as an analog of the Duffin-Schaeffer result about
exponential frames, since it means that Fourier transforms of the
measures $\delta_{x_{j}}$ form a frame on the Fourier transform
side.


\bigskip

\begin{thebibliography}{22}




\bibitem{B64}
A. ~Beurling, {\em Local Harmonic analysis with some applications
to differential operators}, Some Recent Advances in the Basic
Sciences, vol. 1, Belfer Grad. School Sci. Annu. Sci. Conf. Proc.,
A.Gelbart, ed., 1963-1964, 109-125.

\bibitem{BM67}
A. ~Beurling and P. ~Malliavin, {\em On the closure of characters
and the zeros of entire functions}, Acta Math.,\textbf{118},
(1967), 79-95.



\bibitem{BS}
M. Birman and M. Solomyak, {\em Spectral thory of selfadjoint
operators in Hilbert space}, D.Reidel Publishing Co., Dordrecht,
1987.


\bibitem {DS}
R. ~Duffin, A. Schaeffer, {\em A class of nonharmonic Fourier
series}, Trans. AMS, 72, (1952), 341-366.




\bibitem{FG94} H.~Feichtinger and K.~Gr\"ochenig,
{\em Theory and practice of irregular sampling. Wavelets:
mathematics and applications}, 305--363, Stud.\ Adv.\ Math., CRC,
Boca Raton, FL, 1994.



\bibitem{FP} H.~Feichtinger and I.~Pesenson,
{\em Iterative recovery of band limited functions on manifolds},
in Wavelets, Frames and  Operator Theory, Contemp. Math., 345,
AMS, (2004), 137-153.



\bibitem{FP} H.~Feichtinger and I.~Pesenson,
{\em A reconstruction method for  band-limited signals on the
hyperbolic plane}, Sampl. Theory Signal Image Process. 4 (2005),
no. 2, 107--119.


\bibitem {F}
H.~F\"uhr, {\em Abstract Harmonic Analysis of Continuous Wavelet
Transforms,} Lecture Notes in Mathematics, 1863, Springer, 2005.



\bibitem {FG}
H.~F\"uhr and K.~Gr\"ochenig, {\em Sampling theorems on locally
compact groups from oscillation estimates}, preprint, 2005.


\bibitem {Helg1}
S.~Helgason, {\em Differential Geometry and Symmetric
Spaces},Academic, N.Y., 1962.


\bibitem{H}
S.~Helgason, {\em A duality for symmetric spaces with applications
to group representations}, Adv. Math. {\bf 5}, (1970), 1-154.


\bibitem{KPes90}
S.~Krein, I.~Pesenson, {\em Interpolation Spaces and Approximation
on Lie Groups}, The Voronezh State University, Voronezh, 1990,
(Russian).


\bibitem{Lan67}
H.~Landau, {\em Necessary density conditions for sampling and
interpolation of certain entire functions}, Acta. Math.,
\textbf{117}, (1967), 37-52.

\bibitem{LS02}
Y.~Lyubarskii, K.~Seip,{\em Weighted Paley-Wiener spaces}, J.
Amer. Math. Soc. 15(2002),no. 4, 979-1006.


\bibitem{N}
E.~Nelson, {\em Analytic vectors}, Ann. of Math., 70(3), (1959),
572-615.


\bibitem{OS}
J.~Ortega-Cerda, K.~Seip, {\em Fourier frames}, Annals of Math.,
155 (2002), 789-806.

\bibitem{PW34}
R.E.A.C.~Paley and N.~Wiener, {\em Fourier Transforms in the
Complex Domain}, Coll. Publ., \textbf{19},  Providence: Amer.
Math. Soc., (1934).

\bibitem{Pes89}
I.~Pesenson, {\em The Best Approximation in a Representation Space
of a Lie Group}, Dokl. Acad. Nauk USSR, v. 302, No 5, pp.
1055-1059, (1988) (Engl. Transl. in Soviet Math. Dokl., v.38, No
2, pp. 384-388, 1989.)

\bibitem{Pes90}
 I.~Pesenson, {\em  The Bernstein Inequality in the Space of Representation
 of Lie group}, Dokl. Acad. Nauk USSR {\bf 313} (1990), 86--90;
 English transl. in Soviet Math. Dokl. {\bf 42} (1991).

\bibitem{Pes95} I. ~Pesenson, {\em Lagrangian splines,
 Spectral Entire  Functions and Shannon-Whittaker Theorem on Manifolds},
 Temple University Research Report 95-87,  (1995), 1-28.

\bibitem{Pes98a}
I.~Pesenson, {\em  Reconstruction of Paley-Wiener functions on the
Heisenberg group}, Amer. Math. Soc. Transl. (2) Vol. {\bf 184},
(1998), 207- 216.

\bibitem{Pes98b}
I.~Pesenson, {\em Sampling of Paley-Wiener functions on stratified
groups}, J. of Fourier Analysis and Applications {\bf 4} (1998),
269--280.

\bibitem{Pes99}
I.~Pesenson, {\em Reconstruction of band-limited functions in
$L_{2} (R^{d}),$} Proceed. of AMS, Vol.127(12), (1999), 3593-
3600.

\bibitem{Pes00}
 I.~Pesenson, {\em  A sampling theorem on homogeneous manifolds},
 Trans.\ of AMS, Vol. 352(9), (2000), 4257-4270.

\bibitem{Pes01}
I.~Pesenson, {\em Sampling of Band limited vectors}, J. of Fourier
Analysis and Applications {\bf 7}(1), (2001), 93-100 .


\bibitem{Pes04}
I.~Pesenson,  {\em Poincare-type inequalities and reconstruction
of Paley-Wiener functions on manifolds }, J. of Geometric Analysis
{\bf 4}(1), (2004), 101-121.



\bibitem{Pes05}
I.~Pesenson,  {\em Deconvolution of band limited functions on
symmetric spaces}, will appear in the Houston J. of Math.



\bibitem{PP1}
M.~Plancherel, G.~Polya, {\em Fonctions entieres et integrales de
Fourier multiples}, Comment. Math. Helv. 9, (1937), 224-248.


\bibitem{PP2}
M.~Plancherel, G.~Polya, {\em Fonctions entieres et integrales de
Fourier multiples}, Comment. Math. Helv. 10, (1938), 110-163.



\bibitem{T}
A.~Terras,  {\em Harmonic analysis on symmetric spaces and
applications}, Springer-Verlag, 1985.



\bibitem{T2}
H.~Triebel, {\em  Theory of function spaces II,}
  Monographs in Mathematics, 84. Birkhäuser Verlag, Basel, 1992.

\end{thebibliography}
\end{document}